\documentclass[12pt,reqno]{amsart}

\usepackage[latin1]{inputenc}
\usepackage[T1]{fontenc}
\usepackage[british]{babel}
\usepackage{amssymb,graphicx,url,paralist,ifpdf}
\usepackage[margin=30pt,font=small]{caption}

\ifpdf
\fi

\setdefaultenum{(i)}{}{}{}

\theoremstyle{plain}
\newtheorem{theorem}{Theorem}[section]
\newtheorem{proposition}[theorem]{Proposition}
\newtheorem{lemma}[theorem]{Lemma}

\theoremstyle{definition}
\newtheorem{definition}[theorem]{Definition}

\theoremstyle{remark}

\newtheorem{remark}[theorem]{Remark}
\newtheorem*{acknowledgements}{Acknowledgements}

\numberwithin{equation}{section}
\numberwithin{figure}{section}

\newcommand{\Lie}[1]{\operatorname{\textsl{#1}}}
\newcommand{\lie}[1]{\operatorname{\mathfrak{#1}}}
\newcommand{\SU}{\Lie{SU}}

\makeatletter
\newcommand{\abs}{\@ifstar{\absManualScale}{\absAutoScale}}
\newcommand{\absAutoScale}[1]{\left\vert#1\right\vert}
\newcommand{\absManualScale}[2][]{\mathopen #1\lvert #2 \mathclose #1\rvert}

\newcommand{\norm}{\@ifstar{\normManualScale}{\normAutoScale}}
\newcommand{\normAutoScale}[1]{\left\Vert#1\right\Vert}
\newcommand{\normManualScale}[2][]{\mathopen #1\lVert #2 \mathclose #1\rVert}

\newcommand{\inp}{\@ifstar{\inpManualScale}{\inpAutoScale}}
\newcommand{\inpAutoScale}[2]{\left< #1, #2 \right>}
\newcommand{\inpManualScale}[3][]{\mathopen #1\langle #2, #3
\mathclose #1\rangle}
\makeatother

\newcommand{\ii}{\mathbf i}

\newcommand{\tS}{{\mathcal S}}
\newcommand{\CS}{C(\tS)}
\newcommand{\muS}{\mu_\tS}

\newcommand{\Mcut}{M_{\textup{cut}}}

\newcommand{\hkmod}[1]{{#1}_{\textup{mod}}}
\newcommand{\Mmod}{\hkmod{M}}
\newcommand{\Umod}{\hkmod{U}}
\newcommand{\Smod}{\hkmod{\tS}}

\newcommand{\nf}[1]{{#1}_{\textup{non-free}}}
\newcommand{\Mnf}{\nf{M}}

\newcommand{\phiHS}{\phi_{\textup{\textsc{hs}}}}

\newcommand{\bdash}{-\hspace{0pt}} 

\begin{document}
\title{Modifying hyperkähler manifolds with circle symmetry}

\author{Andrew Dancer}
\address[Dancer]{Jesus College\\
Oxford\\
OX1 3DW\\ 
United Kingdom}
\email{dancer@maths.ox.ac.uk}

\author{Andrew Swann}
\address[Swann]{Department of Mathematics and Computer Science\\
University of Southern Denmark\\
Campusvej 55\\
DK-5230 Odense M\\
Denmark}
\email{swann@imada.sdu.dk}

\begin{abstract}
  A construction is introduced for modifying hyperkähler manifolds with
  tri-Hamiltonian circle action, that in favourable situations increases
  the second Betti number by one.  This is based on the symplectic cut
  construction of Lerman.  In $4$ or $8$ dimensions the construction may be
  interpreted as adding a D6-brane.  A number of examples are given and
  generalisations to three-Sasaki and hypersymplectic geometry discussed.
\end{abstract}



\maketitle

\begin{center}
  \begin{minipage}{0.7\linewidth}
    \begin{small}
      \tableofcontents
    \end{small}
  \end{minipage}
\end{center}

\newpage
\section{Introduction}

The symplectic cut construction of Lerman~\cite{Lerman:cuts} has
proved to be a valuable tool in the study of symplectic manifolds with
Hamiltonian circle action.  The purpose of this paper is to
investigate analogues of this construction in quaternionic geometry,
in particular for hyperkähler manifolds.

We define a modification construction for a hyperkähler manifold~$M$ with
tri\bdash Hamiltonian circle action, which involves replacing a level set
of the moment map by its quotient by the circle action.  The relationship
between the complement of the level set in~$M$ and the complement of the
quotient in the new manifold~$\Mmod$ is more complicated than in the
symplectic case.  In the symplectic case one has diffeomorphisms between
corresponding components, but only `half' of the original manifold appears.
For the hyperkähler modification all of~$M$ is involved, but all we can
say, in general, is that there is a third space which is a circle bundle
over each complement.

When we perform the construction, we create a new hyperkähler submanifold
of real codimension four, which is a component of the fixed set of the
circle action in the modified space; this submanifold is in fact the
quotient of a level set mentioned above and is the hyperkähler quotient
of~$M$ at the chosen level.  In physical terms this can be interpreted as
adding a new brane to the space.
 
The lowest dimensional example of this construction produces the
Gibbons-Hawking multi-instanton spaces out of flat space. At each
stage a circle level set is collapsed to a point, but the topology at
long range is also changed. This can be generalised in a
higher-dimensional setting to the toric hyperkähler manifolds
of~\cite{Bielawski-Dancer:toric}, where our construction can be
interpreted as adding a new affine flat to the combinatorial data
associated to the manifold.

We study how the topology changes under our construction.  In
particular, we prove that in the simply-connected case the second
Betti number increases by one.

We conclude by briefly discussing analogous constructions for three-Sasaki
and hypersymplectic spaces.

\begin{acknowledgements}
  The authors thank Brian Steer and Ulrike Tillmann for useful
  discussions, and the Isaac Newton Institute, Cambridge, and the
  Erwin Schrödinger Institute, Vienna, for hospitality during the
  later stages of this work.
\end{acknowledgements}

\section{Symplectic cuts}

We shall first review Lerman's symplectic cut
construction~\cite{Lerman:cuts}.

Let $(M, \omega)$ be symplectic with a Hamiltonian circle action whose
moment map is $\mu\colon M \rightarrow \mathbb R$. We now consider $M
\times \mathbb C$ with the product symplectic structure and the
$S^1$-action
\begin{equation*}
  e^{\ii \theta}  \colon (m, z) \mapsto (e^{\ii \theta} .m,  e^{-\ii \theta}z).
\end{equation*}
whose moment map is
\begin{equation*}
  \Phi \colon (m, z) \mapsto \mu(m) - \abs z^2.
\end{equation*}
The \emph{symplectic cut}~$\Mcut=\Phi^{-1}(\varepsilon)/S^1$ is the
symplectic quotient of~$M\times\mathbb C$ at level~$\varepsilon$.
This fits into the following diagram:
\begin{equation*}
  M  \stackrel{\pi}{\leftarrow} \Phi^{-1}(\varepsilon) 
  \stackrel{p}{\rightarrow} \Mcut.
\end{equation*}
Here $\pi \colon \Phi^{-1}(\varepsilon) \rightarrow M$ is projection
$(m,z) \mapsto m$, and its image is $\{ m \in M : \mu(m) \geqslant
\varepsilon \}$.  The map $p \colon \Phi^{-1}(\varepsilon) \rightarrow
\Mcut$ is just the quotient map for the $S^1$-action.  If the circle
action on~$M$ (or, more generally, on~$\Phi^{-1}(\varepsilon)$), is free
then the fibres of~$p$ are all circles.  However the fibres of~$\pi$
over points in its image are only circles away from the
locus~$\mu^{-1}(\varepsilon)$.  Over~$\mu^{-1}(\varepsilon)$, by contrast,
$\pi$~is injective.

Notice that $\pi$ admits a section $s \colon m \mapsto (m,
+\sqrt{\mu(m) - \varepsilon})$.  This is essentially derived from the
section of the map $\phi_0 \colon z \mapsto \abs z^2$, the moment map
for the action on~$\mathbb C$, which expresses~$\mathbb C$ as a circle
bundle over~$\mathbb R_{\geqslant0}$ with one special (point) fibre
over the origin.  Now $p \circ s$ gives a surjection from $\{m \in M :
\mu(m) \geqslant \varepsilon \}$ onto~$\Mcut$, which is a
diffeomorphism away from~$\mu^{-1}(\varepsilon)$ but has circle fibres
on~$\mu^{-1}(\varepsilon)$.

We can rephrase this as follows.  The set $\Phi^{-1}(\varepsilon)$ is
the disjoint union $\Sigma_1 \cup \Sigma_2$ where
\begin{gather*}
  \Sigma_1 = \{\, (m,z) : \mu(m) > \varepsilon, \  \abs z = +\sqrt{\mu(m)
    -\varepsilon}\, \},\\
  \Sigma_2 = \{\, (m, 0) : \mu(m) = \varepsilon \,\}.
\end{gather*}
On~$\Sigma_1$, each orbit of the circle action contains a unique
$(m,z)$ with $z$ real and positive ($z= +\sqrt{\mu(m) - \varepsilon}$).
Hence $\Sigma_1 /S^1$ may be identified with $\{ m : \mu(m) > \varepsilon
\}$.  On the other hand, $\Sigma_2/S^1$ is just the symplectic
quotient $\mu^{-1}(\varepsilon)/S^1$.

Thus one sees that the symplectic cut $\Phi^{-1}(\varepsilon)/S^1$ may
be viewed as coming from~$\{ m : \mu(m) \geqslant \varepsilon \}$ by
factoring out the circle action on the
boundary~$\mu^{-1}(\varepsilon)$.

\section{Hyperkähler modifications}

To define a hyperkähler manifold we start with a Riemannian
manifold~$(M,g)$ and three compatible complex structures $I$, $J$ and~$K$
such that $IJ=K=-JI$.  Compatibility means that the tensors
$F_I(X,Y)=g(IX,Y)$, etc., are all $2$-forms.  One says that $M$ is
\emph{hyperkähler} if these three $2$-forms are closed.  In particular,
$F_I$, $F_J$ and $F_K$ define symplectic structures on~$M$.

To generalise the symplectic cut construction to hyperkähler geometry,
we replace the extra factor~$\mathbb C$ in the previous section with
$\mathbb H\cong\mathbb R^4$.  The flat metric on~$\mathbb H$ is
hyperkähler with complex structures induced by multiplication by unit
quaternions on the right.  Fixing~$\ii$ we will often write $\mathbb
H=\mathbb C+\mathbf j\mathbb C$ with complex coordinates~$(z,w)$.  This
carries a circle action preserving the geometry given by
\begin{equation}
  \label{eq:action-H}
  e^{\ii \theta}\colon (z,w) \mapsto (e^{\ii \theta}z,e^{-\ii \theta}w).
\end{equation}
This action is \emph{tri\bdash Hamiltonian}, meaning that it is
Hamiltonian with respect to each of the three symplectic structures.
We may combine the three moment maps into the single map taking values
in $\mathbb R^3 = \mathbb R\times\mathbb C$:
\begin{equation}
  \label{eq:mu-H}
  (z,w)\mapsto(\tfrac12(\abs z^2-\abs w^2),\ii zw).
\end{equation}

Now let $M$ be a hyperkähler manifold with a tri-Hamiltonian circle
action.  Consider $M \times \mathbb H$ with action:
\begin{equation*}
  e^{\ii \theta} \colon (m, z,w)
  \mapsto (e^{\ii \theta} m, e^{-\ii \theta}z, e^{\ii \theta}w).
\end{equation*}
The associated moment map $\Phi \colon M \times \mathbb H \rightarrow
\mathbb{R}^3 = \mathbb R \times \mathbb C$ is defined by
\begin{equation}
  \label{eq:Phi}
  \Phi \colon (m, z, w) \mapsto (\mu_{\mathbb R}(m) -\tfrac12 (\abs
  z^2 - \abs w^2),\, \mu_{\mathbb C}(m) - \ii zw),
\end{equation}
where $\mu = (\mu_\mathbb R, \mu_\mathbb C) \colon M \rightarrow
\mathbb R^3 = \mathbb R \times \mathbb C$ is the hyperkähler moment
map for the circle action on~$M$.

\begin{definition}
  The \emph{modification} of the tri\bdash Hamiltonian
  hyperkähler manifold~$M$ at level~$\varepsilon$ is defined to be
  $\Mmod = \Phi^{-1}(\varepsilon)/S^1$, where
  $\varepsilon=(\varepsilon_{\mathbb R}, \varepsilon_{\mathbb C})$ and
  $\Phi$ is as in equation~\eqref{eq:Phi}.
\end{definition}

The $S^1$-action on $M \times \mathbb H$ is free except at
points~$(m,0,0)$ where $m$~is a point of~$\Mnf$, the set of points
of~$M$ with non-trivial stabiliser.  The results
of~\cite{Hitchin-KLR:hK} on the hyperkähler quotient construction,
show that the modification~$\Mmod$ will be a smooth manifold provided
we choose the level~$\varepsilon$ to be outside~$\mu(\Mnf)$.

\begin{definition}
  A \emph{good modification} of~$M$ is the modification~$\Mmod$ at a
  level~$\varepsilon$ lying in~$\mu(M)\setminus\mu(\Mnf)$.  Such
  an~$\varepsilon$ will also be called \emph{good}.
\end{definition}

Any good modification~$\Mmod$ carries a hyperkähler structure obtained by
restricting the Kähler forms of $M\times\mathbb H$ to the level set and
descending to the quotient.  This is a consequence of the general theory of
hyperkähler quotients in~\cite{Hitchin-KLR:hK}.  In addition, in this case
$\Mmod$~will be complete as long as $M$~is complete.

Note that, analogously to the symplectic case, the circle action
\begin{equation}
  \label{eq:action-Mmod}
  (m,z,w) \mapsto (e^{\ii \psi}.m, z, w)
\end{equation}
on $M \times \mathbb H$ preserves the level set of~$\Phi$ and commutes
with our previous $S^1$-action on that set.  So the
action~\eqref{eq:action-Mmod} descends to a tri-Hamiltonian action
on~$\Mmod$.  More specifically, writing $[m,z,w]$ for the point
of~$\Mmod$ represented by~$(m,z,w) \in \Phi^{-1}(\varepsilon)$, the
moment map for the action on~$\Mmod$ is $[m,z,w] \mapsto \mu(m)$.

Summarising the discussion so far we have:

\begin{proposition}
  \label{prop:mod}
  Let $M$ be a hyperkähler manifold with tri\bdash Hamiltonian
  circle action with moment map~$\mu\colon M\to\mathbb R^3$.  Then
  each good modification~$\Mmod$ of~$M$ is a again a smooth
  hyperkähler manifold, of the same dimension as~$M$, with
  tri\bdash Hamiltonian circle action.
\end{proposition}

Note that the circle action on~$\Mmod$ induced
by~\eqref{eq:action-Mmod} is free except when there exist $e^{\ii
  \theta} \in S^1$ and $e^{\ii \psi} \in S^1\setminus\{1\}$ with
\begin{equation*}
  (e^{\ii \psi}.m, z,w) = (e^{\ii \theta}.m, e^{-\ii \theta}z, e^{\ii
    \theta}w).
\end{equation*}
This set is the union of the two sets
\begin{gather*}
  \{ [m,0,0] : \mu(m) = \varepsilon \},\\
  \{ [m,z,w] : m \in \Mnf,\ \Phi(m,z,w) = \varepsilon \}
\end{gather*}
so
\begin{equation*}
  \mu(\nf{(\Mmod)}) = \{ \varepsilon \} \cup \mu(\Mnf).
\end{equation*}
The modification construction may therefore be iterated, provided that we
choose our level set suitably at each stage (and $\mu(M) \setminus
\mu(\Mnf)$ is not a finite set).  In particular, if the circle action
on~$M$ is free then we may iterate the construction provided we change the
level set at each stage.

\begin{remark}
  If one chooses to work in the orbifold category then all that one
  requires is that the circle action be non-trivial and that the
  fixed-point set $M^{S^1}$ does not meet the level
  set~$\mu^{-1}(\varepsilon)$.
\end{remark}

\section{Topology}

Let us now analyse the structure of the modification~$\Mmod$.  Recall
(see~\cite{Bielawski-Dancer:toric} for example) that~\eqref{eq:mu-H} is a
map $\phi \colon \mathbb H \rightarrow \mathbb R^3 = \mathbb R \times
\mathbb C$ such that $\phi^{-1}(0,0)$ is~$(0,0)$ and whose fibre over any
other point is exactly a free orbit of the circle
action~\eqref{eq:action-H}.  This is the hyperkähler analogue of the map
$\phi_0 \colon z \mapsto \abs z^2$ from $\mathbb C$ to~$\mathbb R$
discussed in the symplectic case.  However, there are two important
differences:
\begin{enumerate}
\item \label{item:i}
  the map $\phi$ is onto $\mathbb R^3$, whereas its symplectic
  analogue $\phi_0$ is only onto the closed half\/line $\mathbb
  R_{\geqslant 0}\subset\mathbb R$;
\item \label{item:ii}
  the map $\phi$ does not admit a section, unlike $\phi_0$.
  Indeed on spheres in $\mathbb H$, we see that $\phi$ is the Hopf
  fibration $S^3 \rightarrow S^2$.
\end{enumerate}
Point~(\ref{item:i}) is the reason for our choice of terminology
`modification' instead of `cut' since the whole of~$M$ contributes to
the construction of~$\Mmod$ rather than a proper subset.

We have a diagram
\begin{equation}
  \label{eq:diagram-HK}
  M \stackrel{\pi}{\leftarrow} N \stackrel{p}{\rightarrow} \Mmod.
\end{equation}
Here $N=\Phi^{-1}(\varepsilon)$ and as before $p$~is just the quotient
map for the circle action on~$N$, which has circle fibres if
$S^1$~acts freely on~$N$.  Also $\pi \colon N \rightarrow M$ is just
projection $(m,z,w) \mapsto m$.  It is now onto~$M$ (in contrast to
the symplectic case) because of~(\ref{item:i}).  So we do not remove
half the manifold when we perform the hyperkähler modification.

The fibre of~$\pi$ is a circle except on the set~$\tilde X=\mu^{-1}
(\varepsilon)\times\{(0,0)\}$, where $\pi$~is injective.  Thus under the
modification, the set $X=\mu^{-1}(\varepsilon)$ in~$M$ is being replaced by
its $S^1$-quotient $\hat X=X/S^1$, the hyperkähler quotient of~$M$ at level
$\varepsilon$, in~$\Mmod$.  When $\varepsilon$ is good, this set is a
component of the fixed-point set of the
circle action~\eqref{eq:action-Mmod} on~$\Mmod$, see the discussion after
Proposition~\ref{prop:mod}.  Using the results
of~\cite{Dancer-Swann:geometry-qK}, we thus have:

\begin{proposition}
  Suppose $\Mmod$ is a good modification of~$M$ at level~$\varepsilon$.
  Then $\Mmod$ contains the hyperkähler quotient of~$M$ at
  level~$\varepsilon$ as a hyperkähler submanifold of codimension~$4$.
\end{proposition}

The relation between the complements $M^*=M \setminus X$ and
$\Mmod^*=\Mmod \setminus \hat X$ is more complicated than in the
symplectic case, as (due to~(\ref{item:ii})) the map $\pi$~does not
admit a global section.  What we do have is a space
$N^*=N\setminus\tilde X$ which is a circle bundle over both $M^*$
and~$\Mmod^*$.

\begin{remark}
  If $\mu(M)\ne\mathbb R^3$ then one also has the possibility of
  modifying~$M$ at a level~$\varepsilon\notin\mu(M)$.  In this situation no
  new fixed-points are created, but $\Mmod$ will not be isometric to~$M$
  and again there will be a common circle bundle above both spaces.  This
  situation is excluded by our definition of good modification.
\end{remark}

\subsection{Euler characteristic}
\label{sec:Euler}

As mentioned above the collapsed set~$\hat X$ in~$\Mmod$ will be a
component~$L$ of the fixed-point set of the hyperkähler circle
action~\eqref{eq:action-Mmod} on~$\Mmod$.  It has real codimension~$4$
in~$\Mmod$ and its normal bundle will be a non-trivial circle
representation of quaternionic dimension one.  

Suppose the circle action on~$\Mmod$ has no non-trivial finite
isotropy groups, and that all components of~$L$ are of
codimension~$4$.  Then $\Mmod/S^1$ will be a smooth manifold: the
spheres in the normal bundle are copies of~$S^3$, whose quotient
by~$S^1$ will be~$S^2$.  As this is a boundary it can therefore be
filled in in the quotient, see
\cite{Gukov:Clay,Gukov-S:Spin7,Atiyah-W:G2}.

In the language of string theory such a component of~$L$ is a D6-brane when
the ambient space is $4$- or $8$-dimensional.  Whenever we modify~$M$,
therefore, we are adding in a new brane whose position is determined by the
level~$\varepsilon$.

Note also that, by a standard argument, the Euler characteristic
of~$\Mmod$ will be the sum of the Euler characteristics of the
components of the fixed-point sets of the circle action. Indeed if $L$
has finitely many connected components~$L_i$, for $i=1,\dots,k$, and
$V_i$~are mutually disjoint tubular neighbourhoods of~$L_i$, then as
the circle action on $\Mmod\setminus L$ is free, we see that the
intersection with~$V_i$ has~$\chi=0$, so
\begin{equation*}
  \chi(\Mmod) = \sum_{i=1}^k\chi(V_i) + \chi(\Mmod\setminus L) =
  \sum_{i=1}^k \chi(L_i) = \chi(L). 
\end{equation*}
This can sometimes yield useful information on the topology
of~$\Mmod$.

\subsection{Cohomology}
\label{sec:cohomology}

Our aim is to relate the cohomology of~$M$ and its good
modification~$\Mmod$.  In our cohomology calculations we shall take
coefficients over $\mathbb R$ unless otherwise stated.  We shall also
assume that our manifolds have finite topological type.

We use the notation introduced at the beginning of this section, and in
addition choose a tubular neighbourhood $U$ of~$X=\mu^{-1}(\varepsilon)$
in~$M$.  Noting that $\tilde X=\pi^{-1}(X)$, we write $\tilde
U=\pi^{-1}(U)$ and put $\Umod=p(\tilde U)$.  Again $p(\tilde X)=\hat X$ and
$\Umod$ is a neighbourhood of~$\hat X$ in~$\Mmod$.  As above, a star ${}^*$
on a set will denote the complement of the corresponding set associated
to~$X$, e.g., $\Umod^*=\Umod\setminus\hat X$.

Observe that in the diagram~\eqref{eq:diagram-HK}, $N^*$~is the total
space of circle fibrations over both $M^*$ and~$\Mmod^*$.  We denote
the Euler classes by $e$ and $e^{\prime}$, respectively.  Similarly,
we have that $\tilde{U}^*$~is the total space of circle fibrations
over $U^*$ and~$\Umod^*$.

We note that the normal bundle of~$X$ in~$M$ is trivial.  In fact if
$\xi$ denotes the Killing field for the circle action on $X$, then
$I\xi, J\xi, K\xi$ give an explicit trivialisation.  Hence we may take
$U^*$ to be homotopic to~$S^2 \times X$.

\begin{lemma}
  Suppose $M$ is simply connected.  Then $\Mmod$~is also simply connected
  and
  \begin{equation}
    \label{eq:betti-3}
    b_2(M^*)=b_2(\Mmod^*)=b_2(\Mmod).
  \end{equation}
\end{lemma}

\begin{proof}
  As $X$ has codimension~$3$ in~$M$, we deduce that $M^*$ is also simply
  connected so~$H^1(M^*, \mathbb Z)=0$ and $H^2(M^*,\mathbb Z)$~is
  torsion-free, by the Universal Coefficient Theorem.

  Applying the Gysin sequence to the fibration $S^1 \rightarrow N^*
  \stackrel{\pi} \rightarrow M^*$ we obtain:
  \begin{equation*}
    0 \rightarrow H^1(N^*) \stackrel{\pi_*} \rightarrow H^0 (M^*)
    \stackrel{\wedge e} \rightarrow  H^2(M^*) \stackrel{\pi^*}
    \rightarrow H^2 (N^*)  \stackrel{\pi_*} \rightarrow H^1(M^*)=0
  \end{equation*}
  As $e \neq 0$ we deduce that $H^0(M^*)$ injects into $H^2(M^*)$ so
  $H^1(N^*)=0$ and
  \begin{equation}
    \label{eq:betti-1}
    b_2(M^*) = b_2(N^*) + 1.
  \end{equation}
  The long exact homotopy sequence from $S^1 \rightarrow N^*
  \stackrel{p} \rightarrow M^*$ yields:
  \begin{equation*}
    1 \rightarrow \pi_2(N^*) \rightarrow \pi_2(M^*) \rightarrow {\mathbb Z}
    \rightarrow \pi_1 (N^*) \rightarrow 1; 
  \end{equation*}
  hence $\pi_1(N^*)$ is a homomorphic image of $\mathbb Z$. As
  $H^1(N^*)=0$, we deduce that $\pi_1(N^*)$ is trivial or ${\mathbb Z}/{m
  \mathbb Z}$ for some~$m$.  In the latter case the universal cover~$P$
  of~$N^*$ is a principal circle bundle~$P\to M^*$ such that $P^m=N^*$.
  However, the set~$U\subset M$ is diffeomorphic to the normal bundle
  of~$X$ in~$M$, and over a fibre of this normal bundle the restriction
  of~$N^*$ is homotopic to the Hopf fibration, so $N^*$~admits no such
  $m$th root~$P$.  We conclude that $N^*$ is simply connected.

  The long exact homotopy sequence from $S^1 \rightarrow N^* \rightarrow
  \Mmod^*$ now gives:
  \begin{equation*}
    1 \rightarrow \pi_2(N^*) \rightarrow  \pi_2(\Mmod^*) \rightarrow
    {\mathbb Z} \rightarrow \pi_1(N^*) = 1 \rightarrow \pi_1 (\Mmod^*)
    \rightarrow 1; 
  \end{equation*}
  so $\pi_1 (\Mmod^*)$ is trivial.

  The Gysin sequence from $S^1 \rightarrow N^* \stackrel{p}
  \rightarrow \Mmod^*$ gives:
  \begin{equation*}
    0 \rightarrow H^0 (\Mmod^*) \stackrel{\wedge e^{\prime}} 
    \rightarrow H^2(\Mmod^*) \stackrel{p^*} \rightarrow H^2(N^*)
    \stackrel{p_*} \rightarrow H^1 (\Mmod^*)=0.
  \end{equation*}
  We deduce that $e^{\prime} \neq 0$ and
  \begin{equation}
    \label{eq:betti-2}
    b_2 (\Mmod^*) = b_2(N^*) + 1.
  \end{equation}

  Now $\hat X$ is codimension~$4$ in~$\Mmod$, so
  $\pi_1(\Mmod)=\pi_1(\Mmod^*)=1$ and $\pi_2(\Mmod)=\pi_2(\Mmod^*)$.  By
  Hurewicz $H^2(\Mmod,\mathbb Z)=H^2(\Mmod^*,\mathbb Z)$ and the required
  result follows from \eqref{eq:betti-1} and~\eqref{eq:betti-2}.
\end{proof}

\begin{theorem}
  Suppose $M$ is a hyperkähler manifold with tri\bdash Hamiltonian circle
  action and finite topological type.  Let $\Mmod$ be a good modification
  of~$M$ at a level~$\varepsilon$ such that $X=\mu^{-1}(\varepsilon)$ has
  finite topological type.  If $M$ is simply connected, then
  \begin{equation*}
    b_2(\Mmod) = b_2(M)+1.
  \end{equation*}
\end{theorem}

\begin{proof}
  By the Lemma, we need to compare $b_2(M^*)$ and~$b_2(M)$.

  We have a Thom-Gysin sequence
  \begin{equation}
    \label{eq:TG}
    \dotsb \rightarrow H^i (M) \rightarrow H^i
    (M^*) \rightarrow H^{i-2}(X) \rightarrow H^{i+1}(M) \rightarrow
    \dotsb
  \end{equation}
  This can be obtained by applying Mayer-Vietoris to $M = M^* \cup U$.
  Recalling that $U \simeq X$ and $M^* \cap U = U^* \simeq S^2 \times
  X$, we obtain:
  \begin{multline*}
    \dotsb \rightarrow H^i(M) \rightarrow H^i
    (M^*) \oplus H^i(X)\\ \rightarrow H^{i-2}(X) \oplus H^{i}(X)
     \rightarrow H^{i+1}(M) \rightarrow \dotsb,
  \end{multline*}
  where the $H^{i-2}(X) \oplus H^{i}(X)$ term is the Künneth
  decomposition of $H^i(S^2 \times X)$.  Now if $*$~is a basepoint in
  $S^2$ we can consider $X \stackrel{\kappa} \rightarrow S^2 \times X
  \stackrel{\iota} \rightarrow {\mathbb R}^3 \times X \stackrel{\pi}
  \rightarrow X$, where $\kappa \colon x \mapsto (*, x)$, the map
  $\iota$~is inclusion and $\pi$~is projection. As $\pi \iota \kappa$
  is the identity and $\pi^*$~is an isomorphism on cohomology, we
  deduce $\iota^*$ maps $H^i({\mathbb R}^3 \times X) = H^i(U)=H^i(X)$
  isomorphically onto the $H^i(X)$-part of~$H^i(S^2 \times X)$.  This
  yields the sequence~\eqref{eq:TG}.

  As we are assuming $M$ is simply connected the
  sequence~\eqref{eq:TG} begins
  \begin{equation}
    \label{eq:TG2}
    0 \rightarrow H^2(M) \rightarrow H^2(M^*) \rightarrow
    {\mathbb R}  \rightarrow H^3(M) \rightarrow H^3(M^*) \rightarrow \dotsb
  \end{equation}
  The map from $\mathbb R = H^0(X)\to H^3(M)$ is the Thom isomorphism for
  the normal bundle of~$X$ followed by extension by zero.  So the image of
  $1 \in H^0(X)$ is the closed Poincaré dual~$\eta_X$ of~$X$ in~$H^3(M)$,
  \cite[eqn.~(5.13), p.~51]{Bott-Tu:topology}.  However, the operation of
  taking the Poincaré dual is natural with respect to maps in the sense
  that $\eta_{f^{-1}(S)} = f^* \eta_S$ if $S$ is transverse to~$f$.  This
  gives
  \begin{equation*}
    \eta_X = \eta_{\mu^{-1}(\varepsilon)} = \mu^* \eta_{ \{ \varepsilon \}},
  \end{equation*}
  since $\varepsilon$ is a good and so a regular value of~$\mu$.  Now
  $\eta_{\{\varepsilon\}}$ is zero as it lives in $H^3 (\mathbb R^3)=0$, so
  $\eta_X =0 \in H^3(M)$ and the sequence~\eqref{eq:TG2} splits
  after~$\mathbb R$.  Thus the desired conclusion on the second Betti
  numbers now follows.
\end{proof}

\section{Examples}

\subsection{Multi-instanton metrics}
\label{sec:mi}

Let $M = \mathbb H$ with circle action $q \mapsto e^{\ii t} q$. The
associated moment map $\mu \colon M \rightarrow \mathbb R^3$ has circle
fibres, except for the point fibre over the origin.  Modifying at a
non-zero level, we obtain a new hyperkähler $4$-manifold with a hyperkähler
circle action. One of the level sets of~$\mu$ has been replaced by its
quotient by the circle action, that is, it has been shrunk to a point. So
$\Mmod$ maps to~$\mathbb R^3$ but now has two point fibres. This is the
Eguchi-Hanson space, which has~$b_2=1$.  Repeating the construction gives
the multi-instanton series of Gibbons \& Hawking~\cite{Gibbons-H:multi}.
The topology of these spaces is generated by a chain of $k-1$ two-spheres,
where $k$~is the number of point fibres.  Each time we perform the
construction, therefore, the second Betti number increases by one.

\subsection{Taub-NUT metrics}

Let us take $M$ to be~$\mathbb R^3 \times S^1 = \mathbb H / \mathbb Z$
with the flat hyperkähler structure induced from that on~$\mathbb H$.
The free circle action has a moment map which is just projection onto
the $\mathbb R^3$ factor. Applying the modification construction
shrinks one of the fibres to a point, and we obtain the Taub-NUT space
which is topologically~$\mathbb H$ but has a non-flat metric.
Repeating the construction gives the multi-Taub-NUT series which are
topologically, but not metrically, the same as the multi-instanton
series.

\subsection{Hyperkähler toric manifolds}
\label{sec:HKtoric}

A generalisation of Example~\S\ref{sec:mi} is provided by the toric
hyperkähler manifolds of~\cite{Bielawski-Dancer:toric}. These arise as
hyperkähler quotients of~$\mathbb H^d$ by a sub-torus~$\mathcal N$
of~$\mathbb T^d$.  We choose~$\mathcal N$ by taking its Lie
algebra~$\lie n$ to be the kernel of a surjective linear
map~$\beta\colon \mathbb R^d \to \mathbb R^n$ given by
\begin{equation}
  \label{eq:beta}
  \beta \colon e_k \longmapsto u_k,
\end{equation}
with $u_i\in \mathbb Z^n$ (here $e_1,\dots,e_d$ are the standard basis
vectors for~$\mathbb R^d$).  So $\mathcal N$~is determined by a choice
of~$d$ vectors $u_1, \dots, u_d$ spanning~$\mathbb R^n$.

The hyperkähler quotient~$M$ has dimension~$4n$, and $n= d - \dim
\mathcal N$.  Much of the geometry of~$M$ is determined by the
\emph{flats} $H_k$, $k=1,\dots, d$ which are codimension~$3$ affine
subspaces of~$\mathbb R^{3n}$ defined by equations:
\begin{equation*}
  \{ y  \in \mathbb R^n \otimes \mathbb R^3 : \langle y, u_k \rangle
  = \lambda_k \},
\end{equation*}
where the $\lambda_k = (\lambda_k^1, \lambda_k^2, \lambda_k^3) \in
\mathbb R^3$.  The projections of~$\sum_{k=1}^d \lambda_k^je_k$,
$j=1,2,3$, to~$\lie n^*$ give the choice of level set involved in
defining the hyperkähler quotient~$M$.
 
Now $M$ admits a hyperkähler action of $\mathbb T^n = \mathbb T^d /
{\mathcal N}$.  The moment map $\phi \colon M \rightarrow \mathbb
R^{3n}$ for this action is surjective and its fibres over generic
points are copies of~$\mathbb T^n$. However the fibres are
lower-dimensional tori over points lying in the union of the
flats~$H_k$. More precisely, the $\mathbb T^n$-stabiliser of a point
in~$M$ lying over~$y$ has Lie algebra spanned by those vectors~$u_k$
such that~$y \in H_k$, see~\cite[\S3]{Bielawski-Dancer:toric}.  Also,
$M$~has at worst orbifold singularities if and only if no $n+1$ flats
meet.

The modification~$\Mmod$ may also be viewed as a toric hyperkähler
manifold which is now a quotient of~$\mathbb H^{d+1}$ by a
torus~$\hat{\mathcal N}$ with Lie algebra~$\hat{\lie n}$.  We can
view~$\hat{\lie n}$ as the kernel of a map $\hat{\beta} \colon \mathbb
R^{d+1} \rightarrow \mathbb R^n$, where $\hat{\beta}$~is defined as
in~\eqref{eq:beta}, except that we now introduce an additional vector
$u_{d+1} = \sum_{i=1}^{d} \xi_i u_i$. So $\hat{\lie n}$~is an
extension of~$\lie n$ by the $1$-dimensional algebra~$\langle \xi
\rangle$.  (We can regard $\xi$ as the generator of the circle action
on~$M$ by which we choose to perform the modification).

In terms of the configuration of flats, forming~$\Mmod$ from~$M$
corresponds therefore to adding one new flat~$H_{d+1}$.

Note that the multi-instanton spaces of Example~\S\ref{sec:mi} are
just the $n=1$ case of this discussion. Here the flats~$H_k$ are just
points in~$\mathbb R^3$.

In~\cite[\S6]{Bielawski-Dancer:toric}, the topology of~$M$ was
examined in terms of the combinatorial data of the flats.  We
considered the hyperplanes~$H_k^1$ in~$\mathbb R^n$ defined by
$\langle x, u_k \rangle = \lambda_k^1$.  These divide up $\mathbb R^n$
into a union of polyhedra, and we formed the \emph{bounded polyhedral
  complex}~$\mathcal C$ consisting of the \emph{bounded} polyhedra in
this collection and their faces.

Theorem~6.5 of~\cite{Bielawski-Dancer:toric} showed that $M$~is
simply connected with Poincaré polynomial
\begin{equation}
  \label{eq:Ppoly}
  P_t(M) = \sum_{k=0}^{n} d_k (t^2-1)^k,
\end{equation}
where $d_k$~is the number of $k$-dimensional elements of~$\mathcal C$.
In particular, the odd-dimensional cohomology vanishes.  (Theorem~6.5
is in fact still true if $M$ is allowed to have orbifold
singularities).  It follows from the discussion
of~\S\ref{sec:cohomology} that $b_2(\Mmod) = b_2(M) + 1$.

As an example we can consider the Calabi space $M=T^* \mathbb P^2$,
which corresponds to taking $d=3$, $n=2$, $u_1 = e_1$, $u_2 = e_2$ and
$u_3 = -(e_1 + e_2)$.  The complex $\mathcal C$ therefore consists of
a right triangle and its faces. We have a $\mathbb T^2$-action on~$M$
and can choose any circle subgroup of~$\mathbb T^2$ in forming the
modified space~$\Mmod$.

\begin{figure}[htp]
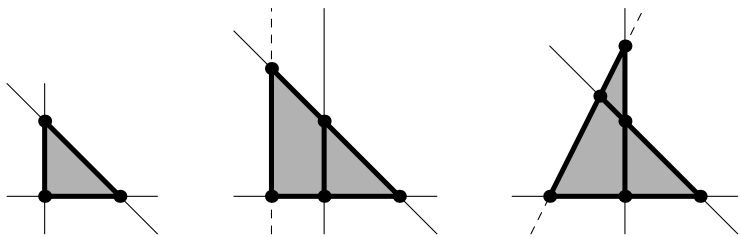

  \centering
  \includegraphics{qcut-pictures.1}
  \qquad
  \includegraphics{qcut-pictures.2}
  \qquad
  \includegraphics{qcut-pictures.3}
  \caption{$T^*\mathbb P^2$ (left) and two of its modifications.}
\label{fig:TP2}
\end{figure}

Forming~$\Mmod$ will involve adding a new hyperplane; depending on the
choice of hyperplane (i.e., the choice of circle action) we can
arrange for $(d_1, d_2)$ of $\Mmod$ to equal~$(6,2)$ or~$(8,3)$, see
Figure~\ref{fig:TP2}.  (For $\Mmod$ to be an orbifold we require that
no three hyperplanes meet).

The Betti numbers of a toric hyperkähler orbifold with $n=2$ are,
from~\eqref{eq:Ppoly}, given by $b_2 = d_1 - 2d_2$ and $b_4 = d_2$.
Note that for $M$ we have $(d_1,d_2)=(3,1)$ so $b_2= b_4=1$, as we
know from the description of $M$ as $T^* \mathbb P^2$.

Both the possibilities for~$\Mmod$ mentioned above yield $b_2=2$ but
$b_4$ can be $2$ or~$3$ depending on the choice of circle action.

The quadratic homogeneity of the moment map involved in forming toric
hyperkähler quotients means that we can apply the arguments for asymptotics
used by Kronheimer~\cite{Kronheimer:ALE} in four dimensions.  Radially
scaling by~$r$ in~$\mathbb H^d$ corresponds to staying instead at the same
distance from~$0$ in~$\mathbb H^d$ and instead scaling~$\lambda$
by~$r^{-2}$.  It follows that the toric hyperkähler quotients are
asymptotically conical in generic directions.  The base of the cone is an
open set in the $3$-Sasaki space corresponding to $\lambda=0$,
cf.~\S\ref{sec:3S}.  However there may be bad directions where this fails,
corresponding to singularities of the $3$-Sasaki quotient
(cf.~Bielawski~\cite{Bielawski:Betti}).

Forming the modification, therefore, preserves the generic
asymptotically conical nature of the manifold, but will in general
alter the base of the cone.

More generally, for general hyperkähler manifolds, we conjecture that
generic asymptotically conical or asymptotically locally conical
(asymptotic to a constant circle over a cone) behaviour will be preserved
under modification, but with the base of the cone being changed, and
possibly with bad non-generic directions being created

\subsection{Gauge theory quotients}

We shall now consider some spaces arising from infinite-dimensional
gauge theoretic constructions.  In contrast to the preceding examples,
these spaces may have odd-dimensional cohomology (in particular
non-zero $b_3$).

Recall from~\cite{Kronheimer:cotangent} that the cotangent
bundle~$T^*G_{\mathbb C}$ of the complexification of a compact Lie
group~$G$ admits a hyperkähler structure.  This is obtained by
identifying~$T^*G_{\mathbb C}$ with the space of $\lie g$-valued
solutions to Nahm's equations smooth on~$[0,1]$, modulo gauge
transformations that are the identity at~$t=0,1$.

These spaces were analysed by the authors
in~\cite{Dancer-Swann:compact-Lie}.  They admit two commuting actions
of~$G$ preserving the hyperkähler structure, corresponding to gauge
transformations that are the identity at one endpoint but not
necessarily at both. It follows
from~\cite[Proposition~4]{Dancer-Swann:compact-Lie} that these actions
are free, as that Proposition constructed equivariant diffeomorphisms
between~$T^* G_{\mathbb C}$ and an open set in~$G \times (\mathfrak
g^*)^3$ containing~$G \times (0,0,0)$, where $G$~acts by left or right
translations on the factor~$G$ in the image.

We therefore have a large supply of free hyperkähler circle actions on~$T^*
G_{\mathbb C}$.  We can perform the modification construction repeatedly
on~$T^* G_{\mathbb C}$, giving new families of hyperkähler manifolds of
dimension $4 \dim G$.  Using a single circle one obtains for example
manifolds with $G\times N(S^1)$-symmetry containing, where $N(S^1)$ is the
normaliser of the chosen circle~$S^1$ in~$G$.  Some of these examples will
contain $T^*(G_{\mathbb C}/\mathbb C^*)$ in the fixed set of the circle
action.  Repeatedly modifying with respect to factors of a maximal
torus~$\mathbb T^n$, one can also obtain families of $G\times \mathbb
T^n$-invariant complete hyperkähler metrics.

Other examples of hyperkähler spaces with free group actions may be
obtained by considering moduli spaces of Nahm data defined on sets of
intervals, where the Nahm matrices are non-singular at some subset of
endpoints.  An example is the hyperkähler $8$-manifold~$M$ with free
$\SU(2)$-action studied in~\cite{Dancer:Nahm-hK,Dancer:hK-family}.
This space is homeomorphic to~$\SU(2) \times \mathbb R^5$.  The
hyperkähler quotients~$\mu^{-1}(\varepsilon)/S^1$ of~$M$ by a circle
subgroup of~$\SU(2)$ all have the homotopy type of the double cover of
the Atiyah-Hitchin manifold, which corresponds to
taking~$\varepsilon=0$.  This non-compact $4$-manifold retracts
onto~$S^2$.  We deduce from \S\ref{sec:Euler} that the corresponding
modification~$\Mmod$ has Euler characteristic equal to~$2$.

\subsection{Toric modifications}

In the symplectic case, Burns, Guillemin \&
Lerman~\cite{Burns-GL:cuts} have generalised the cut construction by
considering symplectic manifolds~$M$ with a~$\mathbb T^n$ action and
moment map~$\mu$.  They then take a toric variety~$X$ of complex
dimension~$n$ with a $\mathbb T^n$-action and associated moment map
$\psi \colon X \rightarrow \Delta$, where $\Delta$ is a polytope
in~$\mathbb R^n$, and consider the product $M \times X$ with the
anti-diagonal $\mathbb T^n$-action and associated moment map $\mu -
\psi$.  The symplectic quotient is then called the symplectic cut
of~$M$ by the toric variety $X$. One can think of the cut as being
obtained by removing the complement of $\mu^{-1} (\Delta +
\varepsilon)$ and factoring out on the boundary the actions of the
torus stabilisers corresponding to faces of~$\Delta$.

We may perform an analogous construction in the hyperkähler case,
taking~$X$ to be one of the \emph{toric hyperkähler manifolds}
mentioned in Example~\S\ref{sec:HKtoric} with quaternionic
dimension~$n$.  As mentioned above the hyperkähler moment map $\phi
\colon X \rightarrow \mathbb R^{3n}$ is surjective, and its fibres are
generically~$\mathbb T^n$ but become lower-dimensional tori on the
intersections~$I_j$ of the flats~$H_k$.  At points of~$X$ in the
pre-image of such intersections, the $\mathbb T^n$-action has
non-trivial torus stabilisers, see~\cite[\S3]{Bielawski-Dancer:toric}.

Forming the generalised hyperkähler modification~$(\Mmod)_{X}$
therefore involves factoring out the sets $\mu^{-1} (I_j +
\varepsilon)$ by the corresponding torus stabilisers.  The complement
of these sets in~$M$ is not necessarily diffeomorphic to the
complement in~$(\Mmod)_X$ of the quotients of these sets; as in the
original construction, all we can say is that there is a third space
which is a $\mathbb T^n$-bundle over each complement.

\section{Modifying three-Sasaki structures}
\label{sec:3S}

A three-Sasaki manifold~$(\tS,g)$ is best described by the property
that the cone~$(\CS=\mathbb R_{>0}\times\tS, dt^2+t^2g)$ is
hyperkähler.  Three-Sasaki manifolds are of dimension~$4n+3$ and
provide many interesting examples of compact Einstein manifolds of
positive scalar curvature; see the survey~\cite{Boyer-G:survey} and
more recent papers by the same authors.

If $\tS$ admits a circle action, then this provides a tri\bdash
Hamiltonian symmetry of the cone~$\CS$.  On $\CS$ there is a unique
choice of hyperkähler moment map~$\mu$ that is homogeneous with
respect to scaling in the $t$~variable.  The restriction~$\muS$
of~$\mu$ to~$\tS$ at~$t=1$ may then be used to define a quotient
construction in the three-Sasaki category provided one reduces at the
level~$0$.

To modify~$\tS$ with respect to the circle action, consider the
hyperkähler modification of $\CS$ at level~$0$.  This is the quotient
of the set
\begin{gather}
  \label{eq:level-S}
  \{ (t,s,z,w) \in \mathbb R_{>0}\times\tS\times\mathbb
  C\times\mathbb C :
  t^2\muS(s) = (\tfrac12(\abs z^2-\abs w^2),\ii zw)
\end{gather}
by the action
\begin{equation}
  \label{eq:S-action}
  e^{\ii\theta}\colon(t,s,z,w) \mapsto (t, e^{\ii\theta}.s,
  e^{-\ii\theta}z, e^{\ii\theta}w).
\end{equation}
Since the defining equations in~\eqref{eq:level-S} are homogeneous
with respect to the scaling $(t,s,z,w)\mapsto(\lambda t,s,\lambda
z,\lambda w)$ the quotient space is a hyperkähler cone of a
three-Sasaki space~$\Smod$ that we call the \emph{modification}
of~$\tS$.  

More concretely, $\Smod$~is the quotient of the set of points
in~\eqref{eq:level-S} satisfying
\begin{equation}
  \label{eq:unit}
  t^2+\abs z^2+\abs w^2=1
\end{equation}
by the action~\eqref{eq:S-action}.  This is a smooth three-Sasaki
manifold of the same dimension as~$\tS$ provided the circle action is
free on $\muS^{-1}(0)\subset\tS$.

\begin{proposition}
  Let $\tS$ be a compact three-Sasaki manifold with circle action that
  is free on~$\muS^{-1}(0)$.  The three-Sasaki modification~$\Smod$ is
  also compact and contains a copy of the three-Sasaki quotient
  of~$\tS$.
\end{proposition}

\begin{proof}
  It is sufficient to show that $t$ is bounded away from zero on the
  intersection of~\eqref{eq:unit} with~\eqref{eq:level-S}.

  As $\tS$ is compact, there is a $K>0$ such that
  $\norm{\muS(s)}\leqslant K$ for all $s\in\tS$.  Since
  $(\tfrac12(\abs z^2-\abs w^2),\ii zw)$ has length $\tfrac12(\abs
  z^2+\abs w^2)$, we see that
  \begin{equation*}
    t^2K \geqslant \norm{t^2\muS(s)} = \tfrac12(1-t^2).
  \end{equation*}
  Thus $t^2\geqslant 1/(1+2K)$, as required.

  The three-Sasaki quotient of~$\tS$ arises as the image of the points
  where $z=0=w$.
\end{proof}

Simple examples of this construction are provided by taking~$\tS$ to
be a round sphere~$S^{4n-1}\subset\mathbb H^n$ with the diagonal
circle action.  The modification is then the three-Sasaki quotient of
the sphere $S^{4n+3}$ by the circle with weights $(1,\dots,1,-1)$,
which is shown to be smooth in~\cite{Boyer-G:survey}.  More generally,
if we consider the circle action on~$S^{4n-1}$ with pairwise co-prime
non-zero weights $(p_1,\dots,p_n)$, then the modification is the
smooth $3$-Sasaki quotient of $S^{4n+3}$ by the circle of weights
$(p_1,\dots,p_n,-1)$.  The cohomology calculations of
Bielawski~\cite{Bielawski:Betti} again show that for toric $3$-Sasaki
orbifolds the modification increases the second Betti number by one.

\section{Hypersymplectic cuts}

Hypersymplectic structures were introduced by
Hitchin~\cite{Hitchin:hypersymplectic} and have subsequently been studied,
for example, in \cite{Hull:actions,Fino-PPS:Kodaira,Dancer-S:toric-hs}.  We
begin with a pseudo\bdash Riemannian manifold~$(M,g)$ and three
endomorphisms $I$, $S$ and~$T$ of the tangent bundle such that $I^2=-1$,
$S^2=+1=T^2$ and $IS=T=-SI$.  These endomorphisms should be compatible
with~$g$ in the sense that the tensors $F_A(X,Y)=g(AX,Y)$, for $A=I,S,T$,
are all $2$-forms; this implies that $g$~has signature~$(2n,2n)$.  One says
that $M$~is \emph{hypersymplectic} if these three $2$-forms are closed.  As
a pseudo-Kähler manifold, $(M,g,I)$ is then Ricci-flat and so neutral
Calabi-Yau, the flat model being~$\mathbb C^{n,n}$.  Also, since $F_I$,
$F_S$ and $F_T$ now define symplectic structures on~$M$, one can discuss
tri\bdash Hamiltonian group actions.  The behaviour of hypersymplectic
manifolds, especially as regards the geometry of their moment maps, is in
some respects intermediate between Kähler and hyperkähler geometries.

Given a hypersymplectic manifold with a circle action and moment map
$\mu \colon M \rightarrow \mathbb R^3$, we can form a hypersymplectic
cut $\Phi^{-1}(\varepsilon_{\mathbb R}, \varepsilon_{\mathbb C})/S^1$ where
\begin{equation*}
  \Phi \colon (m, z, w) \mapsto (\mu_{\mathbb R}(m) - \tfrac12(\abs
  z^2 + \abs w^2),\,
  \mu_{\mathbb C}(m) - \ii z \bar{w})
\end{equation*}
is the moment map for the $S^1$-action
\begin{equation*}
  e^{\ii \theta} \colon (m, z, w) \mapsto (e^{\ii \theta}.m,  e^{-\ii
    \theta}z, e^{-\ii \theta}w)
\end{equation*}
on $M \times \mathbb C^{1,1}$.

Recall from~\cite{Dancer-S:toric-hs} that the moment map $\phiHS \colon
\mathbb C^{1,1} \rightarrow \mathbb R^3$ is given by
\begin{equation*}
  \phiHS \colon  (z,w) \mapsto (\tfrac12 ( \abs z^2 + \abs w^2), \ii z
  \bar{w}) 
\end{equation*}
and descends to a two-to-one map from~$\mathbb C^{1,1}/S^1$ onto $\{(a,b)
\in \mathbb R \times \mathbb C : a \geqslant \abs b \}$, branched over the
boundary $a = \abs b$.  Note the $3$-sphere $\abs z^2 + \abs w^2 = 2a$
in~$\mathbb C^{1,1}$ maps to the disc $\{ b : \abs b \leqslant a \}$.  As
the double cover of the disc branched over the boundary is the $2$-sphere,
our map is essentially the Hopf fibration again.

As usual, we have a diagram
\begin{equation*}
  M \stackrel{\pi}{\leftarrow} 
  \Phi^{-1}(\varepsilon)
  \stackrel{p}{\rightarrow} \Mcut
\end{equation*}
where $p$~is the quotient map for the circle action and $\pi$~is
projection to~$M$.

Now the image of~$\pi$ is the subset of~$M$:
\begin{equation*}
  \pi(M)=\{\, m \in  M  :  \mu_{\mathbb R}(m) - \varepsilon_{\mathbb R}
  \geqslant \abs{\mu_{\mathbb C}(m) - \varepsilon_{\mathbb C}} \,\}.
\end{equation*}
So, as in the symplectic case, but unlike the hyperkähler case, we
are indeed removing part of the hypersymplectic manifold, that lying
over the exterior of the cone in~$\mathbb R^3$.

In some cases, depending on the shape of~$\mu(M)$, and the choice of
cone vertex $(\varepsilon_{\mathbb R}, \varepsilon_{\mathbb C})$ the
hypersymplectic cut may give a compactification of~$\mu(M)$ and even
of~$M$, as in the symplectic case.  An example is given by Example~7.2
of~\cite{Dancer-S:toric-hs}.  However the hypersymplectic structure
may degenerate on a locus within the compactification, as indeed
happens with this specific example.  The same example with a different choice of circle
action on~$M=\mathbb C^{1,1}$ and appropriate choice of level shows one can
alternatively obtain a complete non-degenerate hypersymplectic~$\Mcut$ that
is neither flat nor connected.

Observe also that the fibre of~$\pi$ over a point in~$\pi(M)$ where
inequality is strict is two circles.  Also, $\pi$ is injective over
$\mu^{-1}( \varepsilon_{\mathbb R},\varepsilon_{\mathbb
  C})\times\{(0,0)\}$ and has fibre equal to a single circle over
other points of~$\pi(M)$ where equality holds.

As in the hyperkähler case, the circle fibrations in~$\pi$ will be
non-trivial.

\providecommand{\bysame}{\leavevmode\hbox to3em{\hrulefill}\thinspace}
\providecommand{\MR}{\relax\ifhmode\unskip\space\fi MR }
\providecommand{\MRhref}[2]{%
  \href{http://www.ams.org/mathscinet-getitem?mr=#1}{#2}
}
\providecommand{\href}[2]{#2}

\end{document}